\documentclass{article}
\usepackage{amsmath,amssymb,amsthm}
\usepackage[dvipdfmx]{graphicx}
\begin{document}

\theoremstyle{plain}
\newtheorem{thm}{Theorem}[section]
\newtheorem{lem}[thm]{Lemma}
\newtheorem{cl}[thm]{Claim}

\theoremstyle{definition}
\newtheorem{dfn}[thm]{Definition}

\newcommand{\spl}{\setminus \nobreak \! \nobreak \! \nobreak \setminus}
\newcommand{\blank}{\vspace{0.5\baselineskip}}

\title{THE MINIMAL VOLUME ORIENTABLE HYPERBOLIC 3-MANIFOLD WITH 4 CUSPS}
\author{KEN'ICHI YOSHIDA}
\date{}
\maketitle

\begin{abstract}
We prove that the $8^4_2$ link complement is the minimal volume orientable 
hyperbolic manifold with 4 cusps.
Its volume is twice of the volume $V_8$ of the ideal regular octahedron, 
i.e. $7.32... = 2V_8$. 
The proof relies on Agol's argument used to determine the minimal volume hyperbolic 
3-manifolds with 2 cusps. 
We also need to estimate the volume of a hyperbolic 3-manifold with  
totally geodesic boundary which contains an essential surface with non-separating 
boundary. 
\end{abstract}

\section{Introduction}
\label{section:introduction}

For hyperbolic 3-manifolds, their volumes are known to be topological invariants. 
The structure of the set of the volumes of hyperbolic 3-manifolds is known. 

\begin{thm}[J\o rgensen-Thurston's theorem] 
\emph{(Benedetti-Petronio \cite[corollary E.7.1 and corollary E.7.5]{BP92})}
The set of the volumes of orientable hyperbolic 3-manifolds is a well-ordered set 
of the type $\omega ^\omega$ with respect to the order of $\mathbb{R}$. The volume of an 
orientable hyperbolic 3-manifold with $n$-cusps corresponds to an $n$-fold limit 
ordinal.
\end{thm} 

This theorem gives rise to the problem of determining the minimal volume 
orientable hyperbolic 3-manifolds with $n$ cusps. The answers are known 
in the cases where $0 \leq n \leq 2$.

\begin{itemize}
\item In the case where $n = 0$ (closed manifold), 

Gabai, Meyerhoff and Milley \cite{GMM09} showed that the Weeks manifold has the minimal volume. 
Its volume is 0.94.... 

\item In the case where $n = 1$, 

Cao and Meyerhoff \cite{CM01} showed that the figure-eight knot complement and the manifold 
obtained by the (5,1)-Dehn surgery from the Whitehead link complement 
have the minimal volume. 
Their volume is $2.02... = 2V_3$, where $V_3$ is the volume of the ideal regular 
tetrahedron. 

\item In the case where $n = 2$,

Agol \cite{Ag10} showed that the Whitehead link complement 
and the ($-2$,3,8)-pretzel link complement 
have the minimal volume. Their volume is $3.66... = 
4 \sum_{k=0}^{\infty} \frac{(-1)^{k}}{(2k+1)^{2}} = V_8$, where $V_8$ is the volume 
of the ideal regular octahedron. 
\end{itemize}

In the case where $n \geq 3$, Adams \cite{Ad88} showed that the volume of an $n$-cusped 
hyperbolic 3-manifold is not less than $nV_3$. Agol \cite{Ag10} conjectured the following:

\begin{itemize}
\item In the case where $3 \leq n \leq 10$,

the minimally twisted hyperbolic chain link complement has the minimal volume.

\item In the case where $n \geq 11$,

the ($n-1$)-fold covering of Whitehead link complement has the minimal volume.

\end{itemize}

In this paper, we prove this conjecture in the case where $n = 4$. 

\begin{thm}
\label{thm:main}
The minimal volume orientable hyperbolic 3-manifold with 4 cusps is homeomorphic 
to the $8^4_2$ link complement. Its volume is $7.32... = 2V_8$. 
\end{thm}

\begin{figure}
 \centering 
 \includegraphics[width=12cm,clip]{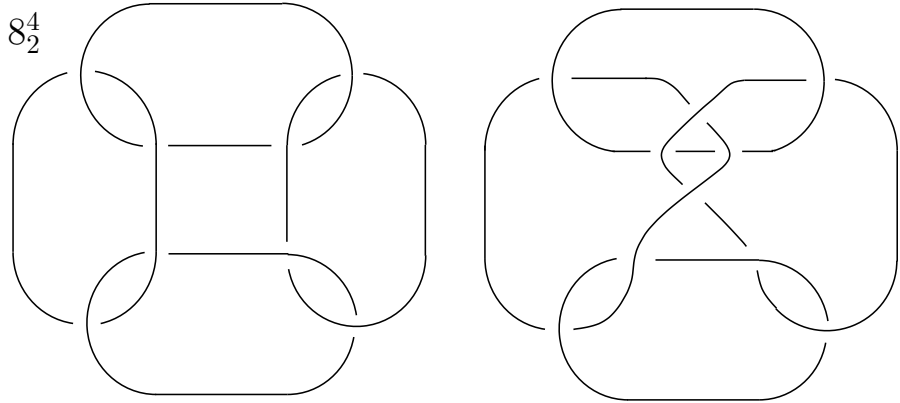} 
 \caption{The $8^4_2$ link and a link whose complement is homeomorphic to 
 that of the $8^4_2$ link} 
\label{fig:links}
\end{figure}

We remark that this link is not the unique one to determine the complement. 
For example, the complement of the link on the right of 
Figure \ref{fig:links} is homeomorphic to the $8^4_2$ link complement. 

We will prove Theorem \ref{thm:main} in Sections \ref{section:estimate} 
and \ref{section:realization}. The proof owes much to Agol \cite{Ag10}.

\section{Review of Agol's argument}
\label{section:review}

In this section, we set up some notation and review the argument used 
by Agol \cite{Ag10} to determine the minimal volume of 2-cusped hyperbolic 3-manifolds. 
We treat compact smooth 3-manifolds with boundary and corners. We only consider 
surfaces in a compact 3-manifold which are properly embedded or 
contained in the boundary. 
Let $I=[0,1]$. 

Let $M$ be a 3-manifold with boundary. For a properly embedded surface $X \subset M$, 
let $M \spl X$ denote the path-metric closure of $M-X$. 
We will say that $X$ is \textit{essential} if $X$ is incompressible and 
$\partial$-incompressible and has no component parallel to the boundary. 
Essential surfaces are not assumed to be connected. 

A finite volume orientable hyperbolic 3-manifold can be the interior of a compact 
3-manifold with the boundary which consists of tori. Its boundary component is called 
a \textit{cusp}. When we say a hyperbolic manifold in what follows, 
it often means this compact manifold. 
We also consider hyperbolic manifolds with totally geodesic boundary. 
In this case there may be \textit{annular cusps} which adjoin the totally geodesic 
boundary. The \textit{double} of a hyperbolic manifold $M$ with totally geodesic boundary 
is the manifold obtained from two copies of $M$ 
by gluing along the totally geodesic boundary. 
Then two annular cusps form one torus cusp in its double. 

\blank

We introduce the notion of pared manifolds. It was defined by Thurston \cite{Th86} to 
characterize a topological property of geometrically finite hyperbolic manifolds. 

\begin{dfn}
(Thurston \cite[Section 7]{Th86}, Morgan \cite[Definition 4.8]{Mo84})
A \textit{pared manifold} is a pair $(M,P)$ such that 

\begin{itemize}
\item 
$M$ is a compact orientable irreducible 3-manifold, 

\item 
$P \subset \partial M$ is a union of annuli and tori which are incompressible in $M$,

\item 
every abelian, noncyclic subgroup of $\pi _1 (M)$ is peripheral with respect to $P$ 
(i.e. conjugate to a subgroup of the fundamental group of a component of $P$), and 

\item 
every map $\phi \colon (S^1 \times I, S^1 \times \partial I) \to (M,P)$ which induces
injective maps on the fundamental groups deforms, as maps of pairs, into $P$. 
\end{itemize}

$P$ is called the \textit{parabolic locus} of the pared manifold $(M,P)$, 
and an annular component of $P$ is called a \textit{pared annulus}. 
We denote by $\partial _0 M$ the surface $\partial M - \mathrm{int} (P)$. 

A pared manifold $(M, P)$ is called \textit{acylindrical} 
if every map $\psi \colon (S^1 \times I, S^1 \times \partial I) \to (M, \partial _0 M)$ 
which induces injective maps on the fundamental groups deforms 
either into $\partial _0 M$ or into $P$. 
\end{dfn}

Since a finite volume orientable hyperbolic 3-manifold is atoroidal, it is 
a pared manifold by letting its parabolic locus be the cusp tori. 
Conversely the following holds: 

\begin{thm}
\label{thm:hyperbolization}
Let $(M,P)$ be an acylindrical Haken pared manifold, and assume that 
$\partial _0 M$ is incompressible. 
We assume that $M$ is not a 3-ball, a $T^2 \times I$ or a solid torus. 
Then $M - P$ admits a finite volume hyperbolic structure with totally geodesic boundary 
$\partial _0 M$. This hyperbolic structure is unique up to isometry. 
\end{thm}

Since the double $DM$ of an acylindrical pared manifold $(M,P)$ is atoroidal, 
$DM$ admits a finite volume hyperbolic structure, 
where $DM$ is obtained from two copies of $M$ by gluing along $\partial _0 M$. 
Then the diffeomorphism swapping the two copies of $M$ can be taken to be an isometry. 
The fixed point set $\partial _0 M$ is totally geodesic \cite[Lemma 2.6]{Le06}. 

\blank

When a hyperbolic manifold is cut along an essential surface, the obtained manifold is 
a pared manifold. 

\begin{lem} 
\label{lem:pared}
\emph{(Agol \cite[Lemma 3.2]{Ag10})}
Let $M$ be a finite volume orientable hyperbolic 3-manifold, and $\partial M$ be 
the parabolic locus $P$ of $M$. Let $X \subset M$ be an essential surface. 
Then $(M \spl X, P \spl \partial X)$ is a pared manifold. 
\end{lem}

\begin{thm}[JSJ decomposition for a pared manifold] 
\label{thm:JSJ}
\emph{(Jaco-Shalen \cite{JS79}, Johannson \cite{Jo79}, Morgan \cite[Section 11]{Mo84})} 
Let $(M,P)$ be a pared manifold such that $\partial _0 M$ is incompressible. 
There is a canonical set of essential annuli $(A, \partial A) \subset 
(M, \partial _0 M)$ called the \textit{characteristic annuli}. 
It is characterized up to isotopy by the property that they are the maximal 
collection of non parallel essential annuli such that every other essential annulus 
$(B, \partial B) \subset (M, \partial _0 M)$ may be relatively isotoped to an 
annulus $(B^{\prime}, \partial B^{\prime}) \subset (M, \partial _0 M)$ 
so that $B^{\prime} \cap A = \emptyset$. 
Then each complementary component $(L, \partial _0 L) \subset 
(M \spl A, \partial _0 M \spl \partial A)$ is one of the following types: 

\begin{enumerate}
\item 
$(T^2 \times I, (T^2 \times I) \cap \partial _0 M)$, where one of the boundary 
components $T^2 \times \partial I$ is a torus component of $P$. 

\item 
$(S^1 \times D^2, (S^1 \times D^2) \cap \partial_0 M)$, which is a solid torus 
with annuli in the boundary. 

\item 
($I$-bundle, $\partial I$-subbundle), which is an $I$-bundle over a surface whose Euler 
characteristic is negative, and the $I$-bundle over the boundary is contained 
in $A \cup P$. 

\item 
$(L, L \cap \partial _0 M)$, where $L$ has no essential annuli whose boundary 
is contained in $L \cap \partial _0 M$. 
\end{enumerate}
\end{thm}

A neighborhood of a torus component of $P$ is either of type 1 or of type 4. 
One of the boundary 
components $T^2 \times \partial I$ of type 1 is a torus component of $P$, and 
the intersection of the other 
boundary component and $\partial _0 M$ is a union of essential annuli in the torus. 
The intersection $(S^1 \times D^2) \cap \partial _0 M$ in a component of type 2 is 
a union of essential annuli in $\partial (S^1 \times D^2)$. 
The union of components of type 3 is called the 
\textit{window}. A component of type 4 is the acylindrical pared manifold 
$(L, L - \partial _0 M)$. The union of the components of type 4 is called 
the \textit{guts} and denoted by $\mathrm{Guts}(M,P)$. 
A torus boundary component of the guts is a torus component of $P$. 

The definition of guts in \cite{Ag10} is a bit different from ours. 
In \cite{Ag10} guts are defined to be the union of types 1, 2 and 4. 
The definition in \cite{AST07} is same as ours, and it is 
appropriate for our purpose. 

Let $M$ be a finite volume orientable hyperbolic manifold. For an essential surface 
$X \subset M$, $(M \spl X, P \spl \partial X)$ is a pared manifold by Lemma \ref{lem:pared}. 
Therefore, 
we can define $\mathrm{Guts}(X) = \mathrm{Guts}(M \spl X, P \spl \partial X)$.  
Then the components of $\mathrm{Guts}(X)$ admit hyperbolic structures 
with geodesic boundary by Theorem \ref{thm:hyperbolization}. 
Hence the volume $\mathrm{vol(Guts}(X))$ is defined. 
This volume is not greater than the volume of $M$. 

\begin{thm}
\label{thm:AST}
\emph{(Agol-Storm-Thurston \cite[Theorem 9.1]{AST07})} 
Let $M$ be a finite volume orientable hyperbolic manifold, and $X \subset M$ be 
an essential surface. Then 
\[
\mathrm{vol}(M) \geq \mathrm{vol(Guts}(X)) \geq \frac{V_8}{2} 
|\chi (\partial \mathrm{Guts}(X))|. 
\]

Moreover, the refinement by Calegari, Freedman and Walker \emph{\cite[Theorem 5.5]{CFW10}}
implies that $M$ is obtained from ideal regular octahedra by gluing along the faces 
when the equality holds. 
\end{thm}

The estimate of $\mathrm{vol(Guts}(X))$ from below in 
Theorem \ref{thm:AST} follows from the following theorem. 

\begin{thm}
\label{thm:Miyamoto}
\emph{(Miyamoto \cite[Theorem 5.2]{Mi94})}
Let $M$ be a hyperbolic manifold with totally geodesic boundary. Then 
$\mathrm{vol}(M) \geq \frac{V_8}{2} |\chi (\partial M)|$. Moreover, $M$ is obtained  
from ideal regular octahedra by gluing along their faces when the equality holds. 
\end{thm}

\begin{lem}
Let $M$ be a finite volume orientable hyperbolic 3-manifold, and $X \subset M$ be 
a non-empty essential surface. Then each component of $\mathrm{Guts}(X)$ has negative 
Euler characteristic. 
\end{lem}

\noindent 
\textit{Proof.}
Since the Euler characteristic of every closed 3-manifold is $0$, 
$\chi (\mathrm{Guts}(X)) \\ = \frac{1}{2} \chi (\partial \mathrm{Guts}(X))$.  
Assume that there is a component $L$ of $\mathrm{Guts}(X)$ such that 
$\chi (L) \geq 0$. Since no component of $\partial \mathrm{Guts}(X)$ is a sphere, 
$\chi (L) = 0$ and $\partial L$ consists of tori. Since $M$ is atoroidal, 
$\partial L \subset \partial M$. This implies $L = M$ by connectedness of $M$. 
This contradicts the fact that $X$ is not empty. 
\hfill $\square$

\blank 

This lemma implies that $\chi (\partial \mathrm{Guts}(X)) \leq -4$ if 
$\mathrm{Guts}(X)$ is not connected. 

\blank

We will use annular compressions to obtain a surface whose guts is not empty.  

\begin{figure}
 \centering 
 \includegraphics[width=12cm,clip]{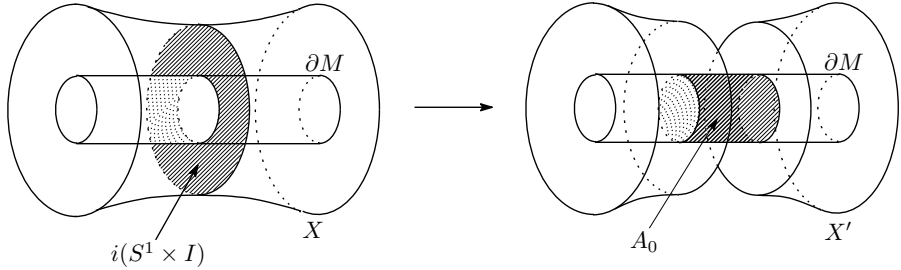} 
 \caption{An annular compression} 
\label{fig:anncomp}
\end{figure}

\begin{dfn}
Let $(X, \partial X) \subset (M, \partial M)$ be an essential surface in a 3-manifold. 
A \textit{compressing annulus} is an embedding $i \colon (S^1 \times I, S^1 \times \{ 0 \}, 
S^1 \times \{ 1 \}) \hookrightarrow (M, X, \partial M)$ such that 

\begin{itemize}
\item 
$i_*$ induces injective maps on $\pi _1$, 

\item 
$i(S^1 \times I) \cap X = i(S^1 \times \{ 0 \} )$, and 

\item 
$i(S^1 \times \{ 0 \} )$ is not isotopic in $X$ to $\partial X$. 

\end{itemize}

An \textit{annular compression} of $(X, \partial X) \subset (M, \partial M)$ is 
the surgery along a compressing annulus $i(S^1 \times I)$. 
Let $U$ be a regular neighborhood of $i(S^1 \times I)$ in $M \spl X$, 
and put $\partial _0 U = \partial U \cap (X \cup \partial M)$ and 
$\partial _1 U = \overline{\partial U - (X \cup \partial M)}$. 
Then the surface $X^{\prime} = (X - \partial _0 U) \cup \partial _1 U$ 
is called the annular compression of $X$. 
If $X$ is essential, $X^{\prime}$ is also essential. 
We will say that $A_0 = \partial U \cap \partial M$ is the annulus in the boundary 
created by the annular compression (Figure \ref{fig:anncomp}). This annulus is not contained 
in the window of $M \spl X^{\prime}$. 
\end{dfn}

\begin{lem}
\label{lem:anncompwindow}
\emph{(\cite[Lemma 3.3]{Ag10})}
Let $M$ be a finite volume orientable hyperbolic manifold. Let $X \subset M$ be 
an essential surface. If $X$ has a compressing annulus, let $X^{\prime}$ be the 
annular compression of $X$. Then the annulus in the 
boundary created by this annular compression is not contained in the window 
of $M \spl X^{\prime}$.
\end{lem}

The following lemma is used in the proof of \cite[Theorem 3.4]{Ag10}.
Lemmas \ref{lem:anncompwindow} and \ref{lem:anncomptorus} imply that 
a torus or an annulus in the boundary is contained 
in the boundary of the gut regions after we perform annular compressions 
as many times as possible. 

\begin{lem}
\label{lem:anncomptorus}
Let $M$ and $X$ be as above. We assume that a $T^2 \times I$ component 
or an $S^1 \times D^2$ component intersects 
a component $T$ of $\partial M$. Then we can perform 
an annular compression for $X$ toward $T$. 
\end{lem}

The following theorem is a result in Culler-Shalen \cite[Theorem 3]{CS84}. 
We will use it to find an essential surface to start 
the proof of Theorem \ref{thm:mainestimate}. 

\begin{thm}
\label{thm:CSsurface}
Let $M$ be a finite volume orientable hyperbolic manifold with $n$ cusps, 
and let $\partial M = T_1 \cup \dots \cup T_n$, where $T_i$ is a torus for 
$1 \leq i \leq n$. Let $k$ be an integer such that $1 \leq k \leq n$. Then there is 
an essential surface $X \subset M$ such that $\partial X \cap T_i \neq \emptyset$ 
for $1 \leq i \leq k$ and $\partial X \cap (T_{k+1} \cup \dots \cup T_n) = \emptyset$. 
\end{thm}

\section{Essential surfaces in 3-manifold with boundary}
\label{section:surface}

In this section we find an essential surface in a hyperbolic 3-manifold with geodesic 
boundary. Using this we will estimate the volume of a hyperbolic 3-manifold with 
geodesic boundary with at least 4 cusps. 
Essential surfaces are found by a homological argument for 3-manifolds, 
and it is not necessary to assume the hyperbolic structure. 

\begin{lem}
\label{lem:homology} 
\emph{(Hatcher \cite[Lemma 3.5]{Ha})}
Let $M$ be a compact orientable 3-manifold. Then the rank of the boundary homomorphism 
$\partial _* \colon H_2 (M, \partial M; \mathbb{Q}) \to H_1 (\partial M; \mathbb {Q})$ 
is half of the dimension of $H_1 (\partial M; \mathbb{Q})$. 
\end{lem}

\begin{lem}
\label{lem:nonsepsurface}
Let $L$ be an orientable hyperbolic 3-manifold with geodesic boundary $S$, with $k$ 
annular cusps $A_1 , \dots , A_k$ and with $n-k$ torus cusps $T_{k+1} , \dots , T_n$, 
where $1 \leq k \leq 3$ and $n \geq 4$. Assume that $\chi (S) = -2$. 
Then there is an essential surface $Y \subset L$ such that 
$Y \cap S = \emptyset$ and $[\partial Y] \neq 0 \in H_1 (\partial L; \mathbb{Z})$. 
\end{lem}

\noindent 
\textit{Proof.} 
The union $S^{\prime} = S \cup A_1 \cup \dots \cup A_k$ is a closed surface of genus 2. 
We note that there are only two types of essential closed curves in $S^{\prime}$, 
one separates $S^{\prime}$ and the other does not. There are no pairs of disjoint
separating curves in $S^{\prime}$. 
 
We can take $k-1$ annuli of $\{ A_1, \dots , A_k \}$ such that the complement 
of them is connected. 
The image of $\partial _* \colon H_2 (L, \partial L; \mathbb{Q}) \to 
H_1 (\partial L; \mathbb{Q})$ is an $(n-k+2)$-dimensional subspace of 
$H_1 (\partial L; \mathbb{Q}) \cong \mathbb{Q}^{2(n-k)+4}$ by Lemma \ref{lem:homology}. 
We consider the subspace $V$ of $H_1 (\partial L; \mathbb{Q})$ spanned by all 
the elements represented by curves in $A_1 , \dots , A_{k-1}, T_{k+1} , \dots , T_n$. 
Since the dimension of $V$ is $2(n-k)+(k-1)$, $V$ intersects 
$\mathrm{Im} (\partial _*)$ in a non-trivial subspace of $H_1 (\partial L; \mathbb{Q})$. 
Hence there exists a non-zero element $z$ in $H_2 (L, \partial L; \mathbb{Q})$ 
such that $\partial _* z \neq 0$ and $z$ belongs to $V$. 
By taking a multiple of $z$, 
there exists a non-zero element $z^{\prime}$ in $H_2 (L, \partial L; \mathbb{Z})$ 
such that $\partial _* z^{\prime} \neq 0$ 
and $\partial _* z^{\prime}$ is represented by a union of 
closed curves in $A_1 , \dots , A_{k-1}, T_{k+1}, \dots , T_n$. 
We can find an essential surface $Y$ representing $z^{\prime}$ such that 
$\partial Y \subset A_1 \cup \dots \cup A_{k-1} \cup T_{k+1} \cup \dots \cup T_n$. 
\hfill $\square$ 

\section{Estimate of volume} 
\label{section:estimate}

Now we are going to estimate the volume of a hyperbolic manifold with geodesic boundary. 
Lemma \ref{lem:nonsepsurface} and Theorem \ref{thm:estimategeodbd} imply that 
the volume of an orientable hyperbolic 3-manifold 
with 4 cusps and with geodesic boundary is not less than $2V_8$. 

\begin{thm}
\label{thm:estimategeodbd} 
Let $L$ be an orientable hyperbolic 3-manifold with geodesic boundary $S$. 
Suppose that there is an essential surface $Y \subset L$ such that 
$Y \cap S = \emptyset$ and $[\partial Y] \neq 0 \in H_1 (\partial L; \mathbb Z)$. 
Then there is an essential surface $Y^{\prime}$ such that 
$\chi (\partial \mathrm{Guts}(L \spl Y^{\prime})) \leq -4$ and 
$\mathrm{vol}(L) \geq 2V_8$. 
\end{thm}

If $\chi (S) \leq -4$, then $\mathrm{vol} (L) \geq 2V_8$ by Theorem \ref{thm:Miyamoto}. 
Hence we may assume that $\chi (S) = -2$. Let $S^{\prime}$ denote the surface which is 
the union of $S$ and the annular cusps of $L$. $\partial L$ consists of 
$S^{\prime}$ and the torus cusps of $L$. 

We will find an essential surface $Y^{\prime}$ such that 
$\chi (\partial \mathrm{Guts}(L \spl Y^{\prime})) \leq -4$. 
Then $\chi (\partial \mathrm{Guts}(DL \spl (DY^{\prime} \cup S)) \leq -8$, 
where $DL$ is the double of $L$ (i.e. the hyperbolic manifold obtained 
from two copies of $L$ by gluing along 
the geodesic boundary $S$) and $DY^{\prime}$ is the union of two copies of 
$Y^{\prime}$ in $DL$. Then Theorem \ref{thm:AST} implies $\mathrm{vol}(DL) \geq 4V_8$, and so 
$\mathrm{vol}(L) \geq 2V_8$. 

We will find a gut component intersecting $S$. For this we need to know 
how a window component intersects $S$. 

\begin{lem}
\label{lem:windowgeodbd}
Let $L$, $S$ and $Y$ be as above. 
Assume that $S$ intersects a component $(J, \partial _0 J)$ of the window of 
$L \spl Y$. 
Then $(J, \partial _0 J)$ is a product $I$-bundle and intersects $S$ 
only on one component of the $\partial I$-bundle. 
\end{lem}

\noindent 
\textit{Proof.} 
Suppose that the base space of $J$ is non-orientable. 
Since $\partial _0 J$ is connected, $\partial _0 J \subset S$. 
We take a simple closed curve $\alpha$ in $J$ such that $\alpha$ is projected to an 
orientation-reversing loop in the base space of $J$. There is a simple closed curve 
$\beta$ in $\partial _0 J$ such that $[\beta ] = [\alpha ]^2 \in \pi _1 (DL) \subset 
\mathrm{Isom}^{+}(\mathbb{H}^3)$. If $\beta$ is homotopic to the boundary of 
$\partial _0 J$, the base space of $J$ is a M\"{o}bius band. It contradicts 
the definition of the window. Hence $[\beta ] \in \pi _1 (S) \subset 
\mathrm{Isom}^{+}(\mathbb{H}^2)$ is hyperbolic element. The simple closed curve $\beta$ 
is homotopic to a simple closed geodesic $\beta ^{\prime}$ in $S$ 
\cite[Theorem 9.6.5]{Ra06}. But the fact that $[\beta ^{\prime}] = [\alpha]^2$ contradicts 
the fact that an element represented by a simple closed geodesic 
in a hyperbolic manifold has no roots \cite[Theorem 9.6.2]{Ra06}. 
Therefore no twisted $I$-bundle intersects $S$. 

Suppose that the base space of $J$ is orientable and both components $Q_0$ and $Q_1$ 
of $\partial _0 J$ are contained in $S$. Since $\chi (Q_0) = \chi (Q_1) < 0$, 
there are (not necessarily simple) closed curves $\gamma _i \subset Q_i$ ($i = 0,1$) 
such that $\gamma _i$ is not homotopic to the boundary of $Q_i$ and 
$\gamma _0$ and $\gamma _1$ are homotopic in $L$. 
Let $\gamma _i ^{\prime}$ be the closed geodesic in $Q_i$ homotopic to $\gamma _i$. 
Since $L$ is totally geodesic, 
the two closed geodesics $\gamma _0 ^{\prime}$ and $\gamma _1 ^{\prime}$ are homotopic 
in $L$. It contradicts the uniqueness of the closed geodesic in a homotopy class. 
Therefore a product $I$-bundle intersects $S$ on at most one side of 
the $\partial I$-bundle. 
\hfill $\square$ 

\blank 

\noindent 
\textit{Proof of Theorem \ref{thm:estimategeodbd}.} 
Let $Y_0$ be an essential surface in $L$ such that $Y_0 \cap S = \emptyset$ and 
$[\partial Y_0] \neq 0 \in H_1 (\partial L; \mathbb Z)$. Moreover let $|\chi (Y_0)|$ 
be minimal among the surfaces satisfying these conditions. 
Since $L$ has no essential sphere, disk, torus or annulus, $\chi (Y_0) < 0$. 
Let $p: L \spl Y_0 \to L$ be the natural projection. 

\blank 

\begin{figure}
 \centering 
 \includegraphics[width=12cm,clip]{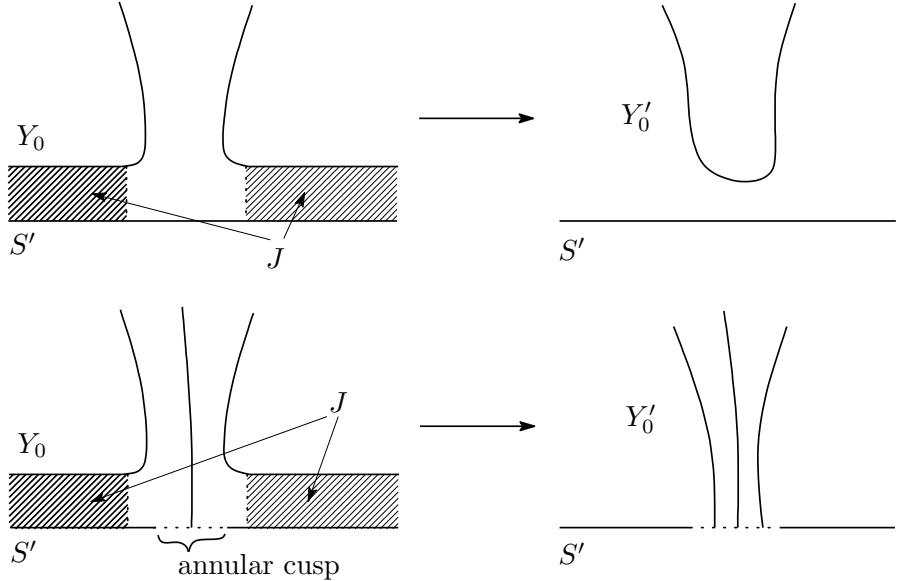} 
 \caption{Constructions in the case where $S$ intersects a component of the window 
 whose Euler characteristic is $-2$} 
\label{fig:construction1}
\end{figure}

(i) First we consider the case where $S$ intersects a component $(J, \partial _0 J)$ 
of the window of $L \spl Y_0$. 
Then $\chi (J)$ is equal to $-1$ or $-2$. 
We will show that $\chi (J) = -1$.  
Assume that $\chi (J) = -2$. 
$S \cap p(J)$ is a 2-punctured torus or a 4-punctured sphere. (If it is a 
closed surface, $Y_0 \cap p(J)$ is a component of $Y_0$ which is parallel 
to $S^{\prime}$. It contradicts that $Y_0$ is essential.) 
Let $Y_0^{\prime}$ be the surface which is 
the union of $Y_0 - (Y_0 \cap p(J))$ and annuli (Figure \ref{fig:construction1}). 
If there is an annulus in $L - Y_0$ whose boundary is two components 
of the frontier of $Y_0 - (Y_0 \cap p(J))$, 
we glue $Y_0 - (Y_0 \cap p(J))$ and this annulus 
(the upper of Figure \ref{fig:construction1}). 
Since $Y_0 \cap p(J)$ is connected, the orientation matches. 
Otherwise, there is an annular cusp which 
is homotopic to the frontier of $Y_0 - (Y_0 \cap p(J))$. 
Then we can glue $Y_0 - (Y_0 \cap p(J))$ and the two annuli, where one of the boundary
components of each annulus is contained in this annular cusp 
(the lower of Figure \ref{fig:construction1}). 
Then $[Y_0^{\prime}, \partial Y_0^{\prime}] = [Y_0, \partial Y_0] \in 
H_2 (L, \partial L; \mathbb{Z})$. We obtain an essential surface from $Y_0^{\prime}$ 
by compressing if necessary. 
Then $|\chi (Y_0^{\prime})| < |\chi (Y_0)|$, 
contradicting the choice of $Y_0$. 
Therefore $\chi (J) = -1$.

\blank 

\begin{figure}
 \centering 
 \includegraphics[width=12cm,clip]{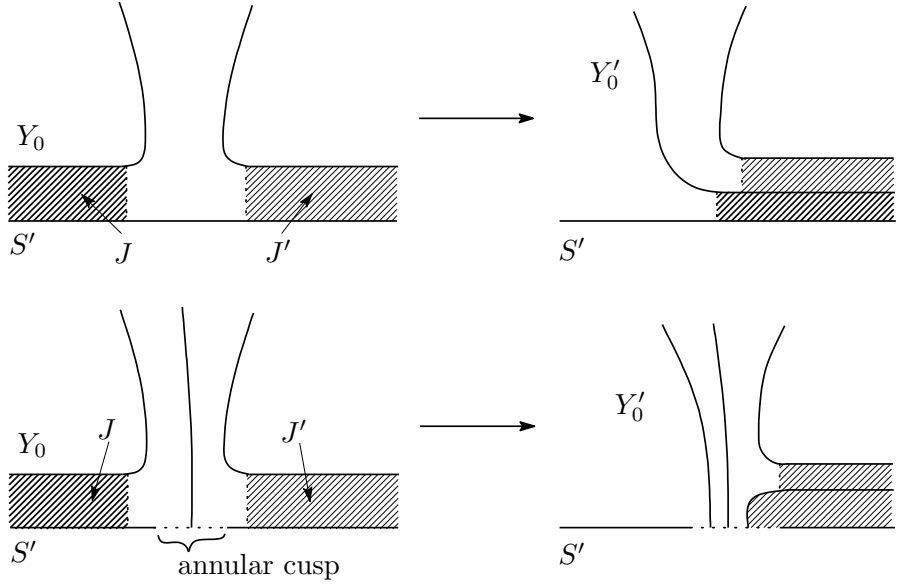} 
 \caption{Constructions in the case where $S$ intersects 2 components of the window 
 whose Euler characteristics are $-1$} 
\label{fig:construction2}
\end{figure}

We will find an essential surface $Y_1$ such that $S$ intersects only one component 
of the window of $L \spl Y_1$. 
If $S$ intersects only one component of the window of $L \spl Y_0$ already, 
put $Y_1 = Y_0$. 
Suppose that $S$ intersects two components $(J, \partial _0 J)$ and  
$(J^{\prime}, \partial _0 J^{\prime})$ of the window of $L \spl Y_0$. 
Let $Y_0^{\prime}$ the surface which is the union of $Y_0 - (Y_0 \cap p(J))$ 
and a surface in $p(J^{\prime})$ (Figure \ref{fig:construction2}). 
Then $[Y_0^{\prime}, \partial Y_0^{\prime}] = [Y_0, \partial Y_0] \in 
H_2 (L, \partial L; \mathbb{Z})$. 
Note that since the orientation may not match, 
we cannot construct a surface as in Figure \ref{fig:construction1}. 
If $Y_0^{\prime}$ is not essential, we obtain an essential surface simpler than $Y_0$ 
by compressing $Y_0^{\prime}$. Since it contradicts the choice of $Y_0$, 
$Y_0^{\prime}$ is essential.

Suppose that $S$ intersects two components of the window of $L \spl Y_0^{\prime}$ again. 
Then one of these two components is contained in $p(J^{\prime})$. 
We can perform the above construction again and remove a part of $Y_0^{\prime}$ 
which is contained in the boundary of the window. 
Since no $I$-bundle can intersect $S$ essentially along both components of the boundary
by Lemma \ref{lem:windowgeodbd}, the part of the obtained surface in $p(J^{\prime})$ 
is not contained in the boundary of the component of the window which intersects $S$ 
and lies on the same side as $p(J)$. 
Hence the above construction can be performed only finitely many times. 

Let $Y_1$ be the essential surface obtained by performing the above construction
as many times as possible. 
The Euler characteristic of the intersection of $S$ and the window of $L \spl Y_1$ 
equals $-1$. 
Therefore the Euler characteristic of the intersection of $S$ and 
$\mathrm{Guts}(L \spl Y_1)$ is equal to $-1$. 
In particular, $\mathrm{Guts}(L \spl Y_1) \neq \emptyset$. 

\blank 

We will find an essential surface $Y_2$ such that 
$\chi (\partial \mathrm{Guts}(L \spl Y_2)) \leq -4$. 
If $\chi (\partial \mathrm{Guts}(L \spl Y_1)) \leq -4$, 
put $Y_2 = Y_1$. 
Suppose that $\chi (\partial \mathrm{Guts}(L \spl Y_1)) = -2$. 
Since the Euler characteristic of $\partial \mathrm{Guts}(L \spl Y_1) - S^{\prime}$ 
is equal to $-1$, it is either a 1-punctured torus or a 3-punctured sphere. 

\begin{figure}
 \centering 
 \includegraphics[width=12cm,clip]{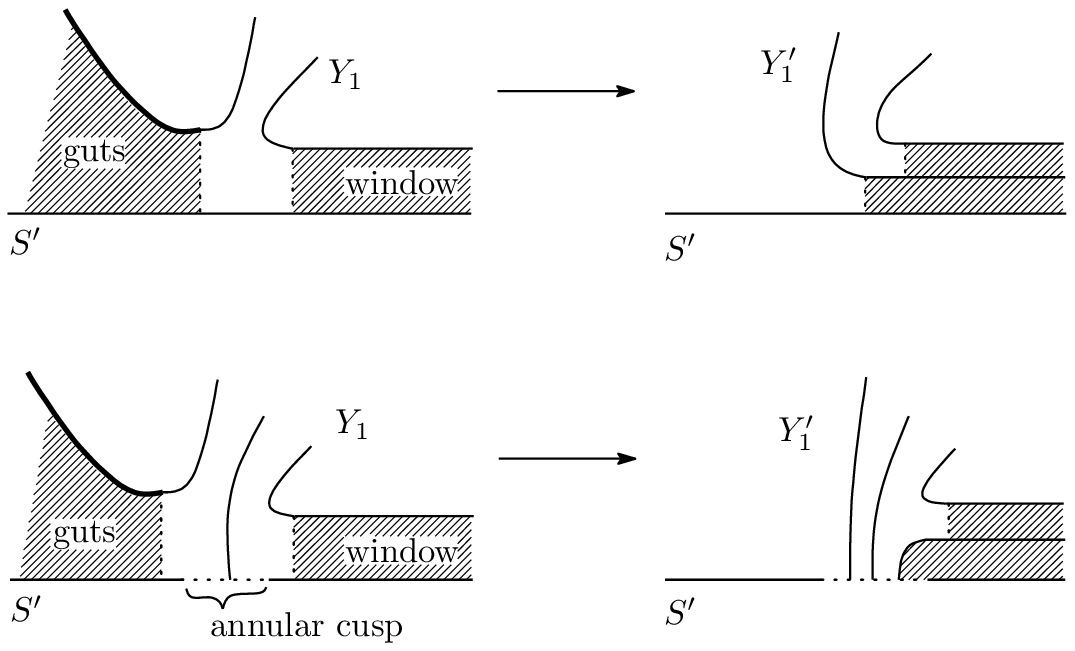} 
 \caption{Constructions in the case where $\partial \mathrm{Guts} (L \spl Y_1) = -2$} 
\label{fig:construction3}
\end{figure}

\begin{figure}
 \centering 
 \includegraphics[width=12cm,clip]{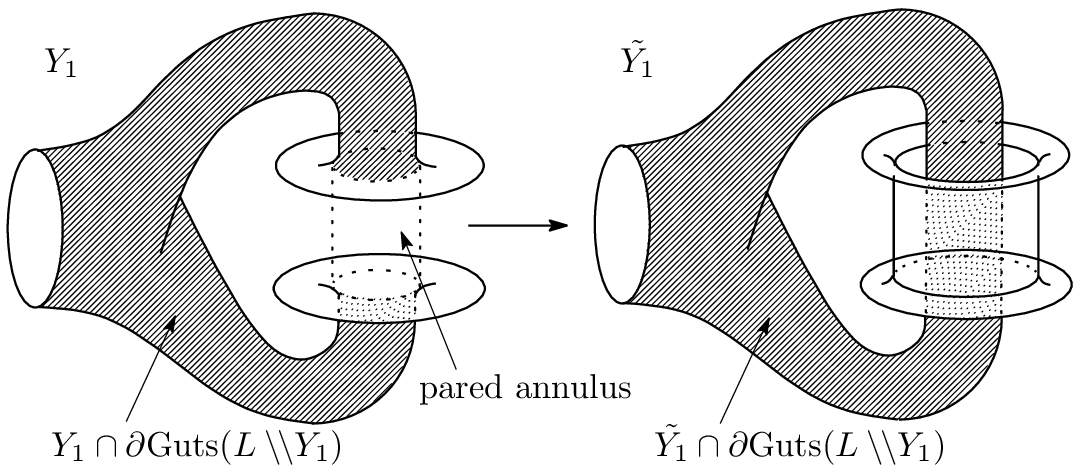} 
 \caption{A construction around a pared annulus in 
 $\partial \mathrm{Guts} (L \spl Y_1) - S^{\prime}$} 
\label{fig:modification}
\end{figure}

Suppose that $\partial \mathrm{Guts}(L \spl Y_1) - S^{\prime}$ is a 1-punctured torus. 
Then 
$\partial \mathrm{Guts}(L \spl Y_1) \\ - S^{\prime}$ 
can contain a pared annulus, 
and $Y_1 \cap \partial \mathrm{Guts}(L \spl Y_1)$ is 
a 1-punctured torus or a 3-punctured sphere. 
If $Y_1 \cap \partial \mathrm{Guts}(L \spl Y_1)$ is a 1-punctured torus, 
let $Y_1^{\prime}$ be the surface which is the union of 
$Y_1 - (Y_1 \cap \partial \mathrm{Guts}(L \spl Y_1))$ and a surface in $p(J^{\prime})$ 
(Figure \ref{fig:construction3}). If $Y_1 \cap \partial \mathrm{Guts}(L \spl Y_1)$ is a 3-punctured sphere,
we obtain the surface $\tilde{Y_1}$ by modifying $Y_1$ 
around the pared annulus in $\partial \mathrm{Guts}(L \spl Y_1) - S^{\prime}$ 
(Figure \ref{fig:modification}). Here $\tilde{Y_1} \cap \partial \mathrm{Guts}(L \spl Y_1)$ is 
a 1-punctured torus. 
Thus we obtain an essential surface $Y_1^{\prime}$ as the union of  
$\tilde{Y_1} - (\tilde{Y_1} \cap \partial \mathrm{Guts}(L \spl Y_1))$ 
and a surface in $p(J^{\prime})$ (Figure \ref{fig:construction3}). 

Suppose that $\partial \mathrm{Guts}(L \spl Y_1) - S^{\prime}$ is a 3-punctured sphere. 
$\partial \mathrm{Guts}(L \spl Y_1) - S^{\prime}$ does not contain a pared annulus. 
Let $Y_1^{\prime}$ be the surface which is the union of 
$Y_1 - (Y_1 \cap \partial \mathrm{Guts}(L \spl Y_1))$ 
and a surface in $p(J^{\prime})$. 

We have obtained a surface $Y_1^{\prime}$ in these ways. 
Then $[Y_1^{\prime}, \partial Y_1^{\prime}] \neq [Y_1, \partial Y_1] 
\in H_2 (L, \partial L; \mathbb{Z})$ in general, but 
$[\partial Y_1^{\prime}] = [\partial Y_1] \neq 0 \in H_1 (\partial L; \mathbb{Z})$. 
Since $|\chi (Y_1)| = |\chi (Y_0)|$, $Y_1$ is essential.

Since $Y_1 \cap \partial \mathrm{Guts} (L \spl Y_1)$
cannot be contained in the window of $L \spl Y_1^{\prime}$, 
$S \cap \partial \mathrm{Guts} (L \spl Y_1)$ is not contained in the window of 
$L \spl Y_1^{\prime}$. 
Hence we can consider that $\mathrm{Guts} (L \spl Y_1^{\prime})$ contains 
$S \cap \partial \mathrm{Guts} (L \spl Y_1)$. 
Since $Y_1$ is essential, 
the added surface in the window is not 
contained in $Y_1^{\prime} \cap \partial \mathrm{Guts} (L \spl Y_1^{\prime})$. 
Hence the above construction can be performed only finitely many times. 

Let $Y_2$ be the essential 
surface obtained by performing the above construction as many times as possible. 
Then $\chi (\partial \mathrm{Guts} (L \spl Y_2))$ is no longer equal to $-2$, 
and $\chi (\partial \mathrm{Guts} (L \spl Y_2)) \leq -4$. 

\blank 

(ii) Suppose that $S$ intersects no component of the window of $L \spl Y_0$. 
Then $\chi (\mathrm{Guts} (L \spl Y_0) \cap S) = -2$. 
Assume that $\chi (\partial \mathrm{Guts} (L \spl Y_0)) = -2$. 
$\partial \mathrm{Guts} (L \spl Y_0)$ is a closed surface which is the union  
of a surface in $S$ and annuli. since $L$ is atoroidal, 
$\partial \mathrm{Guts} (L \spl Y_0)$ 
contains the closed surface $S^{\prime}$. 
Hence $\partial \mathrm{Guts} (L \spl Y_0)$ consists of $S^{\prime}$ and some 
torus cusps of $L$. 
The connectivity of $L$ implies that $L = \mathrm{Guts} (L \spl Y_0)$. 
It contradicts that $Y_0$ is non-empty. Therefore 
$\chi (\partial \mathrm{Guts} (L \spl Y_0)) \\ \leq -4$. 
\hfill $\square$ 

\blank 

We prove the essential part of Theorem \ref{thm:main}. 

\begin{thm}
\label{thm:mainestimate}
Let $M$ be an orientable hyperbolic manifold with 4 cusps. 
Then $\mathrm{vol} (M) \geq 2V_8$. Moreover, if $\mathrm{vol} (M) = 2V_8$, 
$M$ is obtained from two ideal regular octahedra by gluing along the faces. 
\end{thm} 

\noindent 
\textit{Proof.}
It is sufficient to find an essential surface $X \subset M$ such that 
$\chi (\partial \mathrm{Guts} (X)) \\ \leq -4$. Then Theorem \ref{thm:mainestimate} 
follows from Theorem \ref{thm:AST}. 

Let $T_1 , \dots , T_4$ be the cusps of $M$. We take an essential surface $X_0$ 
such that $X_0 \cap T_1 \neq \emptyset$ and $X_0 \cap  T_i = \emptyset 
(2 \leq i \leq 4)$ by Theorem \ref{thm:CSsurface}. 
We perform annular compressions for $X_0$ as many times as possible 
to obtain an essential surface $X_1$. 
When annular compression is 
performed, the number of boundary components of the surface increases 
and its Euler charactersitic 
does not change. Since the Euler characteristic of each component of 
the obtained essential surface is negative, annular compressions can be performed 
only finitely many times. 

We will show that $\mathrm{Guts} (X_1)$ intersects $T_2 , \dots , T_4$. 
Let $k$ be the number of cusps intersecting $X_1 (1 \leq k \leq 4)$. 
Let $T_1 , \dots , T_k$ be the cusps intersecting $X_1$. 
Let $A_2 , \dots , A_k$ be the annuli in 
$T_2 \spl \partial X_1 , \dots , T_k \spl \partial X_1$ created 
by the last annular compressions 
to $T_2 , \dots , T_k$. 
Since there are no compressing annuli any more, Lemma \ref{lem:anncomptorus} implies that 
$A_2 , \dots , A_k$ are not contained in a solid torus component of 
the JSJ decomposition of $M \spl X_1$ and $T_{k+1} , \dots , T_4$ are not contained 
in a $T^2 \times I$ component of it. 
Since compressing annuli to different cusps can be taken disjointly, 
we may change the order of annular compressions to different cusps. 
By Lemma \ref{lem:anncompwindow}, $A_2 , \dots , A_k$ are not contained in the window of $M \spl X_1$. 
Therefore $A_2 , \dots , A_k , T_{k+1} , \dots , T_4 \subset 
\partial \mathrm{Guts} (X_1)$. 

\blank 

If $\mathrm{Guts} (X_1)$ is not connected, 
then $\chi (\partial \mathrm{Guts} (X_1)) \leq -4$ as desired. 
Suppose that $\mathrm{Guts} (X_1)$ is connected. 
Then $A_2 , \dots , A_k , T_{k+1} , \dots , T_4$ are contained in one component $N$ of 
$M \spl X_1$. 
We will find an essential surface $X_2$ such that $\partial \mathrm{Guts} (X_2)$ 
contains at least 4 pared components. 

(i) Suppose that $(T_1 \spl \partial X_1) \cap N \neq \emptyset$. 
If $N = \mathrm{Guts} (X_1)$, let $A_1$ be an annulus which is a component of 
$(T_1 \spl \partial X_1) \cap N$. 
Otherwise let $A_1$ be an separating annulus of the JSJ decomposition intersecting 
$\mathrm{Guts} (X_1)$. 
In either case, $A_1$ is a pared annulus of $\mathrm{Guts} (X_1)$ different 
from $A_2 , \dots , A_k$. Then it is sufficient to put $X_2 = X_1$. 

(ii) Suppose that $(T_1 \spl \partial X_1) \cap N = \emptyset$. 
Let $X_1^{\prime} = X_1 \cap p(N)$, 
where $p \colon  M \spl X_1 \to M$ is the natural projection. 
Then $X_1^{\prime}$ is an essential surface in $M$ 
and $T_1 \cap X_1^{\prime} = \emptyset$. 
$X_1^{\prime}$ is the union of the components of $X_1$ intersecting $N$. 
If we cannot perform an annular compression for $X_1^{\prime}$ to $T_1$, 
$\mathrm{Guts} (X_1^{\prime})$ contains a neighbourhood of $T_1$ 
which is in the complement of $N$. 
Since $\mathrm{Guts} (X_1^{\prime})$ is not connected, 
$\chi (\partial \mathrm{Guts} (X_1^{\prime})) \leq -4$. 
Then it is sufficient to put $X_2 = X_1^{\prime}$. 

If we can perform an annular compression for $X_1^{\prime}$ to $T_1$, 
we obtain $X_2$ by performing annular compressions to $T_1$ as many times as possible. 
Let $A_1$ be the innermost annulus in $T_1$. 
Since $X_1$ is obtained by performing annular compressions as many times as possible, 
there is no compressing annulus for $X_1^{\prime}$ to 
$A_2 , \dots , A_k , T_{k+1} , \dots , T_4$ in $p(N)$. 
Hence there is no compressing annulus for $X_2$ to 
$A_2 , \dots , A_k , T_{k+1} , \dots , T_4$ in $p(N)$. 
Since the surface which is obtained by filling $X_1^{\prime}$ with $A_2 , \dots , A_k$ 
consists of components of a surface in the process of 
the annular compression from $X_0$ to $X_1$, it is essential. 
This implies that $A_1 , \dots , A_k$ are not contained in the window of 
$M \spl X_2$ by Lemma \ref{lem:anncompwindow}. 
Therefore $A_1 , \dots , A_k , T_{k+1} , \dots , T_4 \subset \mathrm{Guts} (X_2)$. 

\blank

\begin{figure}
 \centering 
 \includegraphics[width=10cm,clip]{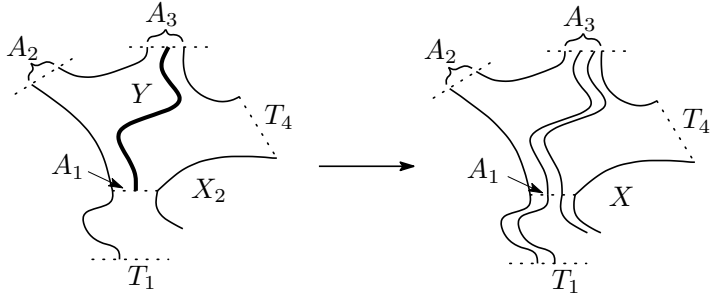} 
 \caption{A construction of an essential surface the boundary of whose guts is 
 no more than $-4$} 
\label{fig:finalconstruction}
\end{figure}

Finally, we will find an essential surface $X$ such that 
$\chi (\partial \mathrm{Guts} (X)) \leq -4$. 
If $k = 4$, the 4 annuli $A_1, \dots , A_4$ are disjoint and not homotopic to each other 
in the non-torus components of $\partial \mathrm{Guts} (X_2)$. 
This implies that $\chi (\partial \mathrm{Guts} (X_2)) \leq -4$. 
Then it is sufficient to put $X = X_2$. 

If $1 \leq k \leq 3$, 
$\mathrm{vol} (\mathrm{Guts} (X_2)) \geq 2V_8$ by Theorem \ref{thm:estimategeodbd}. 
Therefore $\mathrm{vol} (M) \geq 2V_8$ by Theorem \ref{thm:AST}. 
But we need to find $X$ 
in order to prove that $M$ is obtained from 2 octahedra when the equality holds. 
Lemma \ref{lem:nonsepsurface} and Theorem \ref{thm:estimategeodbd} imply that 
there is an essential surface $Y$ in $\mathrm{Guts} (X_2)$ 
such that $\chi (\partial \mathrm{Guts} (\mathrm{Guts} (X_2) \spl Y)) \leq -4$. 
Then $Y$ intersects some of $A_1 , \dots , A_k , T_{k+1} , \dots , T_4$, 
where $A_2 , \dots , A_k , \\ T_{k+1} , \dots , T_4$ 
are contained in $\partial M$. 
If $A_1$ is contained in $\partial M$ or does not intersect $Y$, 
$X_3 \cup Y$ is properly embedded in $M$. 
Since $\mathrm{Guts} (X_2 \cup Y) 
= \mathrm{Guts} (\mathrm{Guts} (X_2) \spl Y)$, 
$\chi (\partial \mathrm{Guts} (X_2 \cup Y)) \leq -4$. 
Then it is sufficient to let $X = X_2 \cup Y$. 
If $A_1$ is contained in the interior of $M$ and intersects $Y$, 
$X_3 \cup Y$ is not properly embedded in $M$. 
Suppose that $A_1 \cap Y$ is the union of $l$ simple closed curves. 
Let $X$ be the union of 2 surfaces parallel to $Y$, 
$X_2 \cap \partial \mathrm{Guts} (X_2)$ and 
$l + 1$ times of $X_2 - \partial \mathrm{Guts} (X_2)$ (Figure \ref{fig:finalconstruction}). 
Since $\mathrm{Guts} (X)$ is homeomorphic to 
$\mathrm{Guts} (\mathrm{Guts} (X_2) \spl Y)$, 
$\chi (\partial \mathrm{Guts} (X)) \leq -4$. 
\hfill $\square$

\section{Realization of hyperbolic manifold}
\label{section:realization}

In this section we will prove that an orientable hyperbolic 3-manifold obtained from
2 ideal regular octahedra by gluing along the faces 
is homeomorphic to the complement of the $8^4_2$ link.
This completes the proof of Theorem \ref{thm:main}. 

Thurston calculated the volume of the complement of the $8^4_2$ link 
in \cite[Ch. 6, Example 6.8.6]{Th78} and it is equal to $2V_8$. 
Moreover, SnapPy \cite{CD} has the list of the orientable hyperbolic 3-manifolds obtained 
from 8 ideal regular tetrahedra by gluing along the faces. 
These imply the uniqueness of the minimal volume 
orientable hyperbolic 3-manifold with 4 cusps, but we prove it here 
by an elementary argument examining the possible ways of gluing along the faces 
of 2 octahedra. 

\blank

The 12 vertices of the 2 octahedra correspond to the 4 cusps of the hyperbolic manifold. 
We look at the number of vertices corresponding to each cusp. 
Since the edge angles of an ideal regular octahedron are right angles, 
4 edges of the 2 octahedra should be glued together. 

\begin{figure}
 \centering 
 \includegraphics[width=12cm,clip]{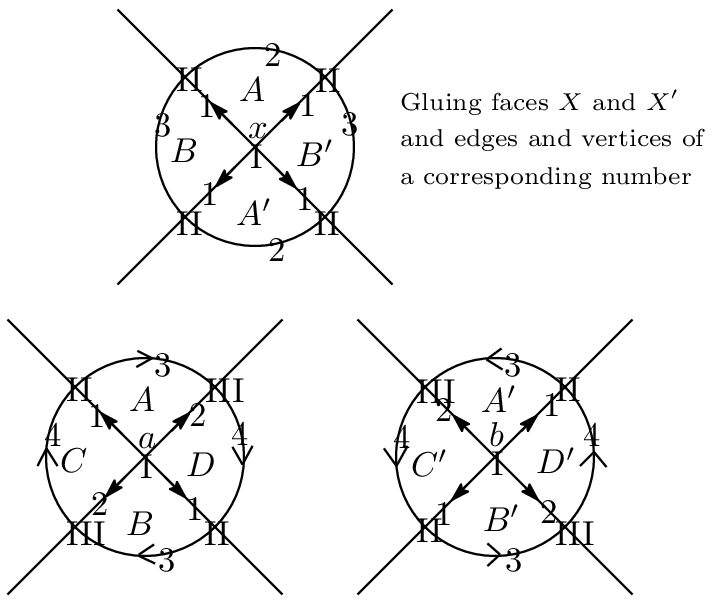} 
 \caption{Gluings of the faces around vertices which are glued together} 
\label{fig:gluingaroundvertices}
\end{figure}

\begin{cl}
\label{cl:gluing}
If there is a cusp consisting of one vertex $x$, the faces around $x$ are glued 
as in the upper part of Figure \ref{fig:gluingaroundvertices}. 
If there is a cusp consisting of 2 vertices $a$ and $b$, the faces around $a$ and $b$ 
are glued as in the lower part of Figure \ref{fig:gluingaroundvertices}. 
\end{cl}

\noindent 
\textit{Proof.} 
If there is a cusp consisting of one vertex $x$, 
the 4 edges around $x$ are glued together, and each face around $x$ is glued with 
the opposite face. 

Suppose there is a cusp consisting of 2 vertices $a$ and $b$. 
Assume that $a$ and $b$ are contained in one octahedron. 
If $a$ and $b$ are adjacent, no edges can be glued with the edge between $a$ and $b$. 
If $b$ is opposite to $a$, we can glue no pairs of faces which are contained in 
different octahedra. This contradicts the connectivity. 
Hence $a$ and $b$ are contained in different octahedra. 

\begin{figure}
 \centering 
 \includegraphics[width=12cm,clip]{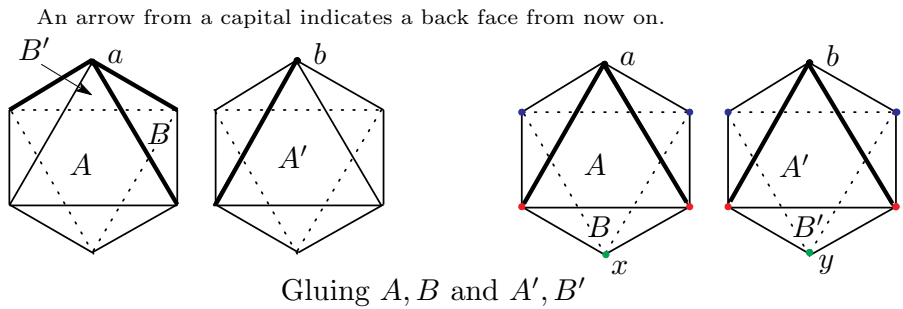} 
 \caption{An impossible example in the case where 2 vertices form a cusp} 
\label{fig:imp2vertices}
\end{figure}

We consider how the 8 edges around $a$ are $b$ are glued. 
Since $a$ and $b$ are glued, the 4 edges around $a$ cannot be glued together. 
If 3 edges around $a$ and one edge around $b$ are glued together, 
2 adjacent faces around $a$ are glued (the left of Figure \ref{fig:imp2vertices}). 
Then the edge between the 2 faces can be glued with no edges. 
Hence 2 edges around $a$ and 2 edges around $b$ are glued together. 
Assume that adjacent edges around $a$ are glued. 
Let $x$ and $y$ be the vertices opposite to $a$ and $b$ respectively. 
If $x$ and $y$ form 2 cusps with themselves, there are 2 edges glued with no other edges. 
Since there are 4 cusps, there is a cusp consisting of $x$ and $y$. 
There are 2 edges glued with no other edges even in this case 
(the right of Figure \ref{fig:imp2vertices}). 
Therefore opposite edges around $a$ are glued and the way of gluing is determined. 
\hfill $\square$

\begin{cl}
\label{cl:no3}
There is no cusp consisting of 3 vertices. 
\end{cl}

\begin{figure}
 \centering 
 \includegraphics[width=10cm,clip]{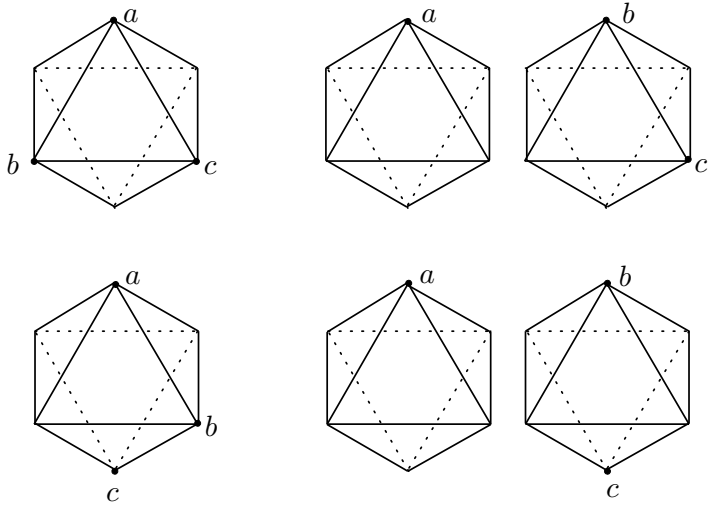} 
 \caption{Positions of 3 vertices} 
\label{fig:positions}
\end{figure}

\noindent 
\textit{Proof.} 
Assume that there is a cusp consisting of 3 vertices $a, b$ and $c$. 
If $a, b$ and $c$ are vertices of one octahedron, 2 positions are possible 
(the left of Figure \ref{fig:positions}). 
If $a, b$ and $c$ are the vertices of one face, 
this face cannot be glued with another face. 
Otherwise, at least one of $a, b$ and $c$ is contained in a face of the octahedron 
containing $a, b$ and $c$. 
This implies that no pair of faces of different octahedra can be glued. 
Hence $a, b$ and $c$ are not contained in one octahedron. 
We assume that $b$ and $c$ are contained in one octahedron without loss of generality. 
Then 2 positions are possible (the right of Figure \ref{fig:positions}). 
If $b$ and $c$ are adjacent, no edges can be glued with the edge between $b$ and $c$. 
Hence $c$ are opposite to $b$. 
Let $x$ be the vertex opposite to $a$. 
Since only the 4 faces around $x$ do not contain $a$, $b$ or $c$, 
the 4 faces cannot be glued with any faces of the other octahedron.

\blank

\begin{figure}
 \centering 
 \includegraphics[width=9cm,clip]{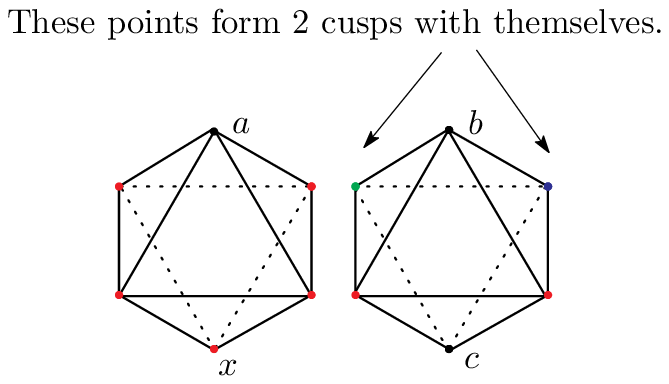} 
 \caption{It is impossible that 5 vertices of an octahedron are glued together.} 
\label{fig:imp5vertices}
\end{figure}

Assume that $x$ does not form a cusp with itself. 
Suppose that adjacent faces around $x$ are glued. 
Then the 5 vertices except $a$ of the octahedron containing $a$ are glued together. 
There are 2 vertices $y$ and $z$ which form 2 cusps with themselves on the octahedron 
containing $b$ and $c$. 
The 4 vertices around $y$ are glued together by Claim \ref{cl:gluing} 
(Figure \ref{fig:imp5vertices}). 
This contradicts that $b$ is glued only with $a$ and $c$. 
Hence opposite faces around $x$ are glued.

\begin{figure}
 \centering 
 \includegraphics[width=8cm,clip]{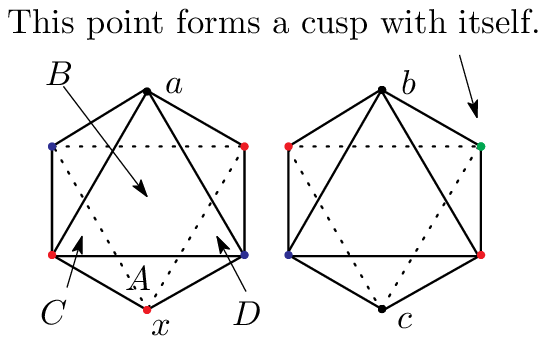} 
 \caption{It is impossible that opposite faces are glued twistedly.} 
\label{fig:imptwist}
\end{figure}

Suppose that opposite faces $A$ and $B$ around $x$ are glued twistedly, i.e. 
the 2 vertices corresponding with $x$ in $A$ and $B$ are not glued. 
Then 2 opposite vertices on the octahedron containing $a$ are glued with $x$. 
Since the 5 vertices except $a$ of the octahedron containing $a$ 
cannot be glued together, 
the other faces $C$ and $D$ around $x$ are glued twistedly. 
The 4 faces around $a$ are glued with faces of the other octahedron 
because of the correspondence of the vertices and the fact that adjacent faces 
around $a$ cannot be glued. 
Hence there is a vertex which forms a cusp with itself on the octahedron containing $b$ and $c$ 
(Figure \ref{fig:imptwist}). This contradicts that $b$ is glued only with $a$ and $c$.

\blank

Hence $x$ forms a cusp with itself. 
Since 4 edges are glued together, 
the faces around $a$ are glued with faces of the other octahedron. 
At least 3 vertices of the octahedron containing $b$ and $c$ are glued 
with the 4 vertices except $a$ and $x$. 
Since we must obtain 4 cusps, there is a vertex which form a cusp with itself. 
It is contradiction. 
\hfill $\square$ 

\blank 

\begin{figure}
 \centering 
 \includegraphics[width=13cm,clip]{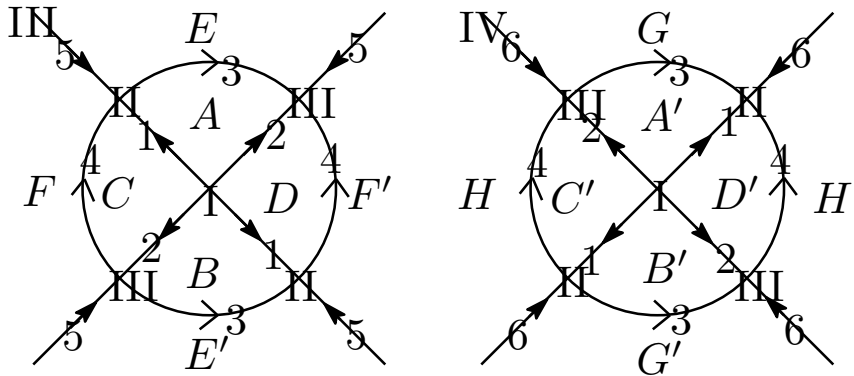} 
 \caption{The way of gluing of the octahedra (i)} 
\label{fig:gluingway1}
\end{figure}

Claim \ref{cl:no3} implies that there is a cusp consisting of one or 2 vertices. 
Suppose that there is a cusp consisting of one vertex $x$. 
The 4 vertices around $x$ are glued together. 
The 4 faces $A, B, C$ and $D$ around the vertex $a$ opposite to $x$ are glued with faces 
of the other octahedron. 
$a$ is glued with only one vertex $b$ because of Claim \ref{cl:no3} and the fact that 
7 vertices are glued. 
Since the 8 vertices around $a$ and $b$ are glued together, the vertex $y$ opposite to 
$b$ forms a cusp with itself. 
The numbers of the vertices corresponding to the cusps are $1,1,2$ and $8$. 
By Claim \ref{cl:gluing} the way of gluing is determined as in Figure \ref{fig:gluingway1} (i).

\begin{figure}
 \centering 
 \includegraphics[width=13cm,clip]{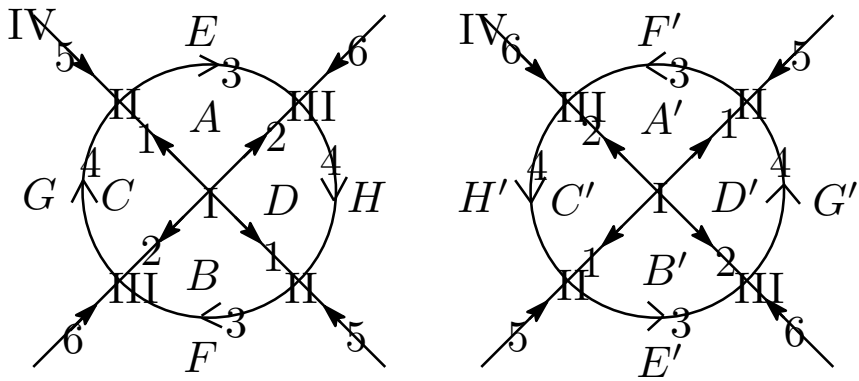} 
 \caption{The way of gluing of the octahedra (ii)} 
\label{fig:gluingway2}
\end{figure}

Suppose that there is no cusp consisting of one vertex. 
Then there is a cusp consisting of 2 vertices $a$ and $b$. 
$A, A^{\prime}, B, B^{\prime}, C, C^{\prime}, D$ and $D^{\prime}$ around $a$ and $b$ 
are glued as Figure \ref{fig:gluingaroundvertices}. 
Since no cusp consists of one vertex, 
the 2 vertices $x$ and $y$ opposite to $a$ and $b$ respectively are glued together. 
The numbers of the vertices corresponding to the cusps are $2,2,4$ and $4$. 
The face $E$ adjacent to $A$ is glued with 
the face $E^{\prime}$ adjacent to $B^{\prime}$ 
because of the correspondence of the vertices and edges. 
The way of gluing is determined as in Figure \ref{fig:gluingway2} (ii).

\begin{figure}
 \centering 
 \includegraphics[width=14cm,clip]{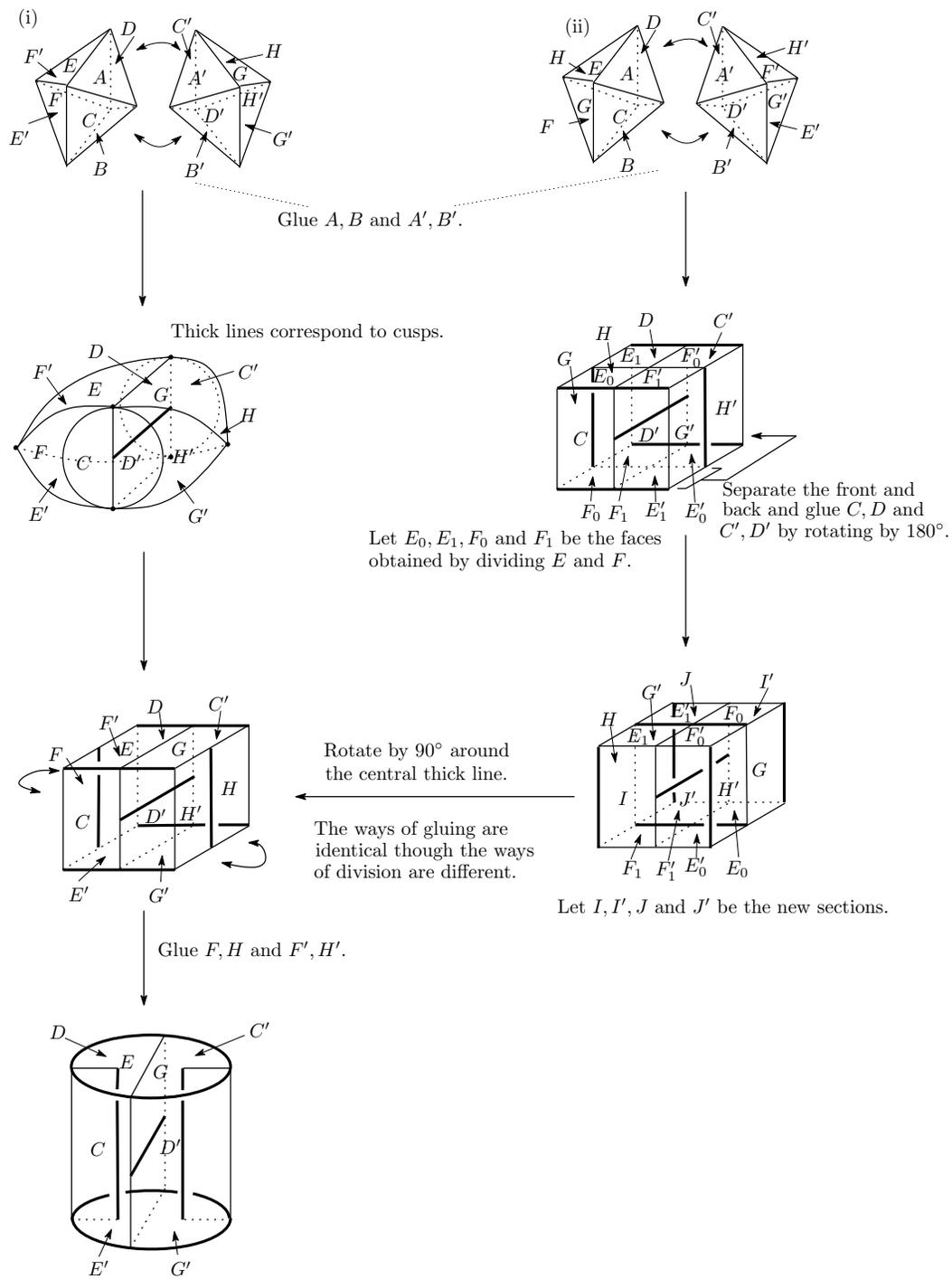} 
 \caption{Gluing of the octahedra I} 
\label{fig:gluingI}
\end{figure}

\begin{figure}
 \centering 
 \includegraphics[width=15cm,clip]{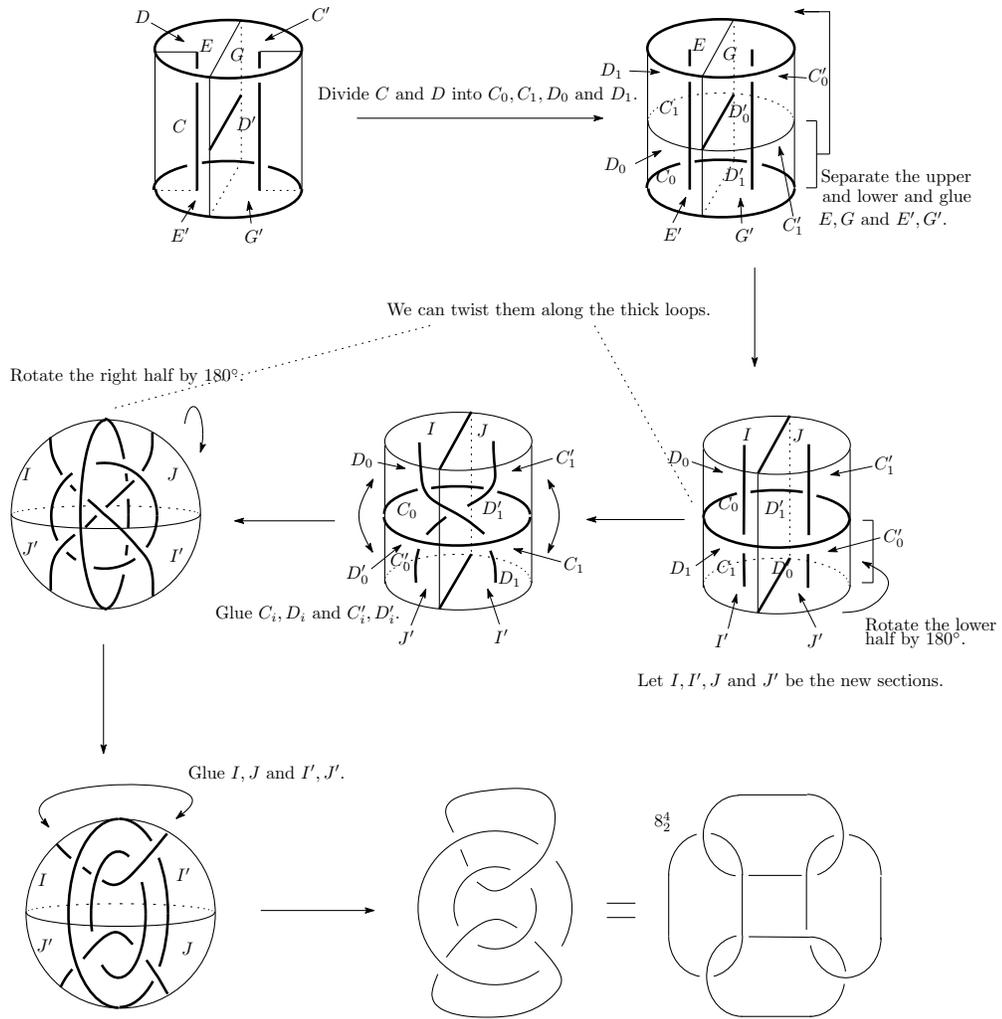} 
 \caption{Gluing of the octahedra II} 
\label{fig:gluingII}
\end{figure}

Both cases of (i) and (ii) give homeomorphic spaces by Figure \ref{fig:gluingI} 
and they are the $8^4_2$ link complements by Figure \ref{fig:gluingII}. 

\clearpage

\clearpage

\clearpage 

\section*{Acknowledgements} 

The author would like to express his gratitude to Takashi Tsuboi 
for helpful guidance and advice. 
This research is supported by Global COE Program 
``New Development in Mathematics'' from the Ministry of Education, Culture, 
Sports, Science and Technology (Mext) of Japan.

\textsc{Graduate School of Mathematical Sciences, 
University of Tokyo, 3-8-1 Komaba, 
Meguro-ku, Tokyo 153-8914, Japan.} 

\textit{E-mail address}: \texttt{kyoshida@ms.u-tokyo.ac.jp}

\end{document}